\documentclass[a4paper,12pt]{article}
\usepackage{amsmath,amssymb}
\usepackage{amscd}
\title{Descent theory for open varieties}
\author{David Harari and Alexei N. Skorobogatov}
\date{\today}
\DeclareFontFamily{U}{wncy}{}
\DeclareFontShape{U}{wncy}{m}{n}{%
   <5>wncyr5%
   <6>wncyr6%
   <7>wncyr7%
   <8>wncyr8%
   <9>wncyr9%
   <10>wncyr10%
   <11>wncyr10%
   <12>wncyr6%
   <14>wncyr7%
   <17>wncyr8%
   <20>wncyr10%
   <25>wncyr10}{}
\DeclareMathAlphabet{\cyrille}{U}{wncy}{m}{n}
\def\Sha{\cyrille X}
\def\Be{\cyrille B}
\def\Ga{\Gamma}
\def\Si{\Sigma}
\def \R{{\bf R}}
\def \A{{\bf A}}
\def \P{{\bf P}}
\def \T{{\cal T}}
\def \K{{\cal K}}
\def \D{{\cal D}}

\def \Pic{{\rm Pic}}
\def \lra{{\longrightarrow}}

\newtheorem{theo}{Theorem}[section]
\newtheorem{prop}[theo]{Proposition}
\newtheorem{lem}[theo]{Lemma}
\newtheorem{cor}[theo]{Corollary}
\newtheorem{defi}[theo]{Definition}

\def \kbar {{\bar k}}

\def \dem {\paragraph{Proof}}

\def \rem {\paragraph{Remark}}
\def \rems {\paragraph{Remarks}}

\def \Romannumeral #1 {\expandafter\uppercase\expandafter {\romannumeral #1} }
\def \br {{\rm{Br}}}

\def \pic {{\rm {Pic}}}
\def \div {{\rm{Div}}}
\def \Div {{\rm{Div}}}

\def \O {{\cal O}}
\def \X {{\cal X}}
\def \Y {{\cal Y}}
\def \S {{\cal S}}
\def \spec {{\rm{Spec\,}}} 
\def \Spec {{\rm{Spec}}} 
\def \Coker {{\rm{Coker\,}}}
\def \Gal {{\rm{Gal}}}

\def\ov{\overline}

\def \Hom {{\rm {Hom}}}
\def \Ext {{\rm {Ext}}}

\def \Z {{\bf Z}}
\def \Q {{\bf Q}}

\def \G {{\bf G}}

\def \hyp{{\mathbb H}}
\def\smallsquare{\vbox{\hrule\hbox{\vrule height 1 ex\kern 1 ex\vrule}\hrule}}
\def\enddem{\hfill \smallsquare\vskip 3mm}

\begin{document}

\maketitle

\begin{abstract}
\noindent We extend the descent theory of Colliot-Th\'el\`ene 
and Sansuc to
arbitrary smooth algebraic varieties by removing the condition
that every invertible regular function is constant. This
links the Brauer--Manin obstruction for integral points on 
arithmetic schemes to the obstructions defined by torsors
under groups of multiplicative type.
\end{abstract}


\medskip

Let $X$ be a smooth and geometrically integral variety 
over a number field $k$ with points everywhere locally.
Descent theory of Colliot-Th\'el\`ene and 
Sansuc \cite {desc}, \cite{skobook}
describes arithmetic properties of $X$ in terms of 
$X$-torsors under $k$-groups of multiplicative type.
It interprets
the {\it Brauer--Manin obstruction} to the existence of 
a rational point (or to weak approximation) on $X$
in terms of the obstructions defined by torsors.

Let $\kbar$ be an algebraic closure of $k$.
Because of its first applications the descent theory was stated in
\cite {desc} for proper varieties
that become rational over $\kbar$; in this case it is enough
to consider torsors under tori. It was pointed out in \cite{S99}
that the theory works more generally under the sole
assumption that the group $\kbar[X]^*$ of 
invertible regular functions on $\ov X:=X \times_k \kbar$ is
the group of constants $\kbar^*$. This assumption is satisfied
when $X$ is proper, but it often fails for
complements to reducible divisors in smooth projective varieties;
it also fails for many homogeneous spaces of algebraic groups.
In the general case of an arbitrary smooth and geometrically
integral variety 
Colliot-Th\'el\`ene and Xu Fei have recently introduced a
Brauer--Manin obstruction to the existence of integral points
\cite[Sect. 1]{ctfei}. 
Descent obstructions to the existence of integral points
were briefly considered by Kresch and Tschinkel in \cite{KT}, Remark 3, 
see also Section 5.3 of \cite{ctwitt}.
In the particular case of an open subset of $\P^1_k$ 
a variant of the main theorem of descent linking
the two kinds of obstructions
has recently turned up in connection with
an old conjecture of Skolem, see \cite[Thm. 1]{dhvoloch}.

The goal of this paper is to 
extend the theory of descent to the general case of a
smooth and geometrically integral variety.
It turns out that the main results are almost entirely the same.
The methods, however, must be completely overhauled. As it
frequently happens, one needs to systematically
consider Galois hypercohomology of complexes instead
of Galois cohomology of individual Galois modules. 
For principal homogeneous spaces of algebraic 
groups this approach has already been 
used in \cite{boduke}, \cite{dhsz2} and \cite{dem2}. But even in the `classical'
case $\kbar[X]^*=\kbar^*$ working with 
derived categories and hypercohomology of complexes streamlines 
the proof of a key result of descent theory
(\cite{desc}, Prop. 3.3.2 and Lemme 3.3.3, \cite{skobook}, Thm. 6.1.2 (a))
by avoiding delicate explicit computations with cocycles (see our
Theorem~\ref{cuppt} and its proof).

\smallskip

Let us now describe the contents of the paper.
Let $S$ be a $k$-group of multiplicative type, that is,
a commutative algebraic group whose connected component of the identity
is an algebraic torus.
In Section \ref{one} we define the {\it extended type} of an $X$-torsor
under $S$, an invariant
that classifies $X$-torsors up to twists by a $k$-torsor.
When $\kbar[X]^*\not=\kbar^*$ the extended type defines a stronger 
equivalence relation on $H^1(X,S)$ 
than the classical type introduced by Colliot-Th\'el\`ene
and Sansuc in \cite{desc}.

Let $\T$ be an $X$-torsor under $S$.
In Section \ref{two} we show that if $U\subset X$ is 
an open set such that the classical type of the torsor
$\T_U\to U$ is zero, then $\T_U$ is canonically isomorphic to
the fibred product $Z\times_Y U$, where $Z$ and $Y=Z/S$ are 
$k$-torsors under groups of multiplicative type, and $U\to Y$ is a 
certain canonical morphism (Theorem \ref{th1}).
This description follows the ideas of Colliot-Th\'el\`ene 
and Sansuc \cite{desc} who used similar constructions to describe
$\T_U$ by explicit equations. Our goal was to obtain a functorial
description, so our results are not immediately related to theirs.
Corollary \ref{co1} describes the restriction of torsors of given 
extended type to sufficiently small open subsets. 

In Section \ref{three} we prove the main results of our generalised
descent theory. The proof of Theorem~\ref{cuppt} 
relies on the previous work of T. Szamuely and 
the first named author \cite{dhsza, dhsz2}, in particular, on
their version of the Poitou--Tate duality for tori, which was later 
extended by C. Demarche to the groups of 
multiplicative type \cite{dempoitou}.

In Section \ref{four}
we prove statements about the existence of integral points and strong 
approximation. As an application we give a short proof of a
result by Colliot-Th\'el\`ene and Xu Fei, generalised by
C. Demarche, see Theorem \ref{miau}.

\section{The extended type of a torsor} \label{one}

Let $Z$ be an integral regular Noetherian 
scheme, and let $p:X\to Z$ be a faithfully flat
morphism of finite type. 
Let $\D(Z)$ be the derived category
of bounded complexes of fppf or \'etale sheaves on $Z$.
For an object ${\cal C}$ of $\D(Z)$, the 
hypercohomology groups $\hyp^i(Z,{\cal C})$ will be
denoted simply by $H^i(Z,{\cal C})$.
Notation such as $\Hom_Z(A,B)$ or $\Ext_Z ^i(A,B)$ will be
understood in the category of sheaves on $Z$, or in $\D(Z)$.
The same conventions apply when $Z$ is replaced by $X$.

Consider the truncated object
$\tau_{\leq 1} \R p_* \G_{m,X}$ in $\D(Z)$. 
Its shift by 1, which has trivial cohomology outside 
the degrees $-1$ and $0$, is denoted by
$$KD(X)=(\tau_{\leq 1} \R p_* \G_{m,X})[1].$$
There is a canonical morphism
$i:\G_{m,Z}\to \tau_{\leq 1} \R p_* \G_{m,X}$, and we define
$$KD'(X)=\Coker(i)[1],$$
so that we have an exact triangle
\begin{equation} \label{triangle}
\G_{m,Z}[1] \lra KD(X)\stackrel{v}{\lra}  KD'(X) 
\stackrel{w}{\lra} \G_{m,Z}[2].
\end{equation}

A group scheme of finite type over $Z$
is called a $Z$-group of multiplicative type if
locally on $Z$ it is isomorphic to a group subscheme of $\G_{m,Z}^n$.
By \cite[\Romannumeral9, Prop. 2.1]{sga3}
such a group is affine and faithfully flat over $Z$.
If $S$ is a group of multiplicative type or a 
finite flat group scheme over $Z$, we denote by 
$\widehat S$ the {\it Cartier dual} of $S$. This is 
the group scheme over $Z$ which represents 
the fppf sheaf ${\cal H}om_Z(S,\G_{m,Z})$, 
see \cite[\Romannumeral10, Cor. 5.9]{sga3} when $S$
is of multiplicative type, and \cite[Ch. 14]{mum} when
$S$ is finite flat.

The following proposition is a generalisation of the
fundamental exact sequence of Colliot-Th\'el\`ene 
and Sansuc (see \cite{skobook}, Thm.~2.3.6 and Cor.~2.3.9). 

\begin{prop} \label{extype}
Let $S$ be a $Z$-group scheme. Assume that one of the two following 
properties is satisfied:

{\rm (a)} $S$ is of multiplicative type;

{\rm (b)} $S$ is finite and flat, and if $2$ is a residual characteristic
of $Z$, then the $2$-primary torsion subgroup $S \{2 \}$ 
is of multiplicative type (equivalently, the Cartier dual 
$\widehat{{S \{2 \}}}$ is smooth over $Z$).  

Then there is an exact sequence
\begin{equation} \label{newtype}
H^1(Z,S) \to H^1(X,S) \stackrel{\chi}{\to} 
\Hom_Z(\widehat S, KD'(X)) \stackrel{\partial}{\to} H^2(Z,S) \to H^2(X,S).
\end{equation}
\end{prop}
To simplify notation, here and elsewhere
we write $H^n(X,S)$ for the fppf  
cohomology group $H^n(X, p^*S)$. If
$S$ is smooth, then the fppf topology can be replaced by \'etale topology.

\paragraph{Proof of Proposition~\ref{extype}} 
We apply the 
functor $\Hom_k(\widehat S,.)$ to the exact triangle (\ref{triangle}).
To identify the terms of the resulting long exact sequence
we use the following well known fact: for any scheme $X/Z$,
any $Z$-group $S$ of multiplicative type and any
$n \geq 0$ we have
$$H^n(X,S)=\Ext_X ^n (p^* \widehat S,\G_{m,X}),$$
see \cite{desc}, Prop.~1.4.1, or \cite{skobook}, Lemma~2.3.7.
Let us recall the argument for the convenience of the reader.
One proves first that ${\cal E}xt^n_X(p^* \widehat S,\G_{m,X})=0$
for any $n\geq 1$, and then the local-to-global spectral
sequence
$$H^m(X,{\cal E}xt^n_X(p^* \widehat S,\G_{m,X}))\Rightarrow
\Ext^{m+n}_X(p^* \widehat S,\G_{m,X})$$
completely degenerates, giving the desired isomorphism.

In case (b) the same argument works
for $1\leq n \leq 3$: indeed, for $\ell\not=2$ we have 
${\cal E}xt^n_X(p^* \widehat S\{ \ell \},\G_{m,X})=0$ 
by the main result of \cite{breen}.
(Note that the case $\ell=2$ is 
exceptional: for example, ${\rm Ext}^2_K (\alpha_2,\G_{m,K}) \neq 0$ 
if $K$ is a 
separably closed field of characteristic $2$, see \cite{breen2}.)

\smallskip

The functor $\R \Hom_X(p^* \widehat S,.)$ from $\D(X)$ 
to the derived category
of abelian groups $\D({\rm Ab})$ is the composition 
of the functors $\R p_*(.):\D(X)\to \D(Z)$ and 
$\R \Hom_Z(\widehat S,.):\D(Z)\to \D({\rm Ab})$.
This formally entails a canonical isomorphism
$$\Ext_X ^n (p^* \widehat S,\G_{m,X})=
R^n \Hom_Z(\widehat S,\R p_* \G_{m,X}).$$ 
In particular, we have 
$$H^n(Z,S)=\Ext^n_Z(\widehat S,\G_{m,Z})=\Hom_Z(\widehat S,\G_{m,Z}[n]).$$

Truncation produces an exact triangle
$$\tau_{\leq 1} \R p_* \G_{m,X} \to 
\R p_* \G_{m,X} \to \tau_{\geq 2} \R p_* \G_{m,X} \to 
(\tau_{\leq 1} \R p_* \G_{m,X}) [1],$$
and here $\tau_{\geq 2} \R p_* \G_{m,X}$ is acyclic in degrees 
$0$ and $1$. We deduce canonical isomorphisms
$$R^1 \Hom_Z (\widehat S,\tau_{\leq 1} \R p_* \G_{m,X})=
R^1 \Hom_Z (\widehat S,\R p_* \G_{m,X})=H^1(X,S),$$ 
and an injection of 
$R^2 \Hom_Z (\widehat S,\tau_{\leq 1} \R p_* \G_{m,X})$ into 
$R^2 \Hom_Z (\widehat S,\R p_* \G_{m,X})=H^2(X,S)$. 
Now (\ref{newtype}) is obtained by applying 
$\Hom_Z(\widehat S,.)$ to (\ref{triangle}).
\enddem

\rems

1. Let $k$ be a field of characteristic zero
with algebraic closure $\kbar$ and Galois group $\Ga=\Gal(\kbar/k)$.
In the case when $X$ is smooth over $k$, $KD(X)$ was introduced in 
\cite{dhsz2} as the following complex of 
$\Ga$-modules in degrees $-1$ and $0$:
$$[\kbar(X)^* \to \div (\ov X)].$$
Here $\kbar(X)$ is the function field of $\ov X=X \times_k \kbar$, 
and $\div (\ov X)$ is the group of divisors on $\ov X$ (see
\cite{bovh}, Lemma~2.3 and Remark~2.6). In this case $KD'(X)$ 
is quasi-isomorphic to the complex of $\Ga$-modules
$$[\kbar(X)^*/\kbar^* \to \div (\ov X)].$$
Up to shift, $KD'(X)$ was independently introduced
by Borovoi and van Hamel in \cite{bovh}: in their notation we have
$KD'(X)=U\pic(\ov X)[1]$.
Furthermore, if $X^c$ is a smooth compactification of $X$, and 
$\div_{\infty}(\ov X^c)$ is the group of divisors of $\ov X^c$ 
supported on $\ov X - \ov X^c$, then, by \cite{dhsz2}, Lemma 2.2, 
$KD'(X)$ is quasi-isomorphic to the complex 
$$[\div_{\infty} (\ov X^c) \to \pic (\ov X^c)].$$

2. In the relative case, 
when $X$ is smooth over $Z$, 
our $KD(X)$ and $KD'(X)$ coincide with 
analogous objects defined in 
\cite{dhsz2}, Remark 2.4 (2). See Appendix A to the present paper
for the proof of this fact. 

\smallskip

3. In the relative case, when $X$ is proper over $Z$ with geometrically 
integral fibres, $KD'(X)$ identifies with the sheaf 
$R^1p_* \G_{m,X}$, the relative Picard functor. 
When $X$ is also assumed projective over $Z$, the relative Picard functor
is representable by a $Z$-scheme, 
separated and locally of finite type,
see \cite{blr}, Ch. 8, Thm. 1 on p. 210. 

\bigskip

Let $\D(k)$ be the bounded derived category of the category
of continuous discrete $\Ga$-modules.

\begin{defi} 
{\rm Let $X$ be a smooth and geometrically integral variety
over $k$. Let $Y$ be an $X$-torsor under a $k$-group of multiplicative
type $S$, and let $[Y]$ be its class in $H^1(X,S)$.
We shall say that the 
morphism $\chi([Y]): \widehat S \to KD'(X)$ in the derived category 
$\D(k)$ is the {\it extended type} of the torsor $Y\to X$.}
\end{defi}

\rems  1. There is a canonical morphism from (\ref{newtype}) to the
sequence of
Colliot-Th\'el\`ene and Sansuc (\cite{skobook}, Thm.~2.3.6):
\begin{equation} \label{classtype}
\Ext^1_k(\widehat S,\kbar[X]^*) \to  H^1(X,S) \to 
\Hom_k(\widehat S,\pic (\ov X)) \to \Ext^2_k(\widehat S,\kbar[X]^*) \to 
 H^2(X,S)
\end{equation}
Indeed, (\ref{classtype}) is obtained by applying the functor 
$\Hom_k(\widehat S,.)$ to the exact triangle
\begin{equation} \label{triangle2}
\kbar[X]^*[1] \to KD(X) \to \pic (\ov X) \to \kbar[X]^*[2]
\end{equation}
(cf. \cite{skobook}, p. 26), and there is an obvious canonical 
morphism from (\ref{triangle}) to (\ref{triangle2}).
Recall that if $Y\to X$ is a torsor under $S$, then the image of
the class $[Y]\in H^1(X,S)$ in $\Hom_k(\widehat S,\pic (\ov X))$
is called the {\it type\/} of $Y\to X$. We see 
that the notion of extended
type defines a stronger equivalence relation on $H^1(X,S)$ than
the notion of type.
For example two torsors have the same extended type if 
and only if their classes in $H^1(X,S)$ coincide up to a `constant element'.

\smallskip

2. If we assume further that $\kbar[X]^*=\kbar^*$ (e.g. $X$ 
proper), then 
$KD'(X)$ is quasi-isomorphic to $[0 \to \pic (\ov X)]$, and
the exact sequence
(\ref{newtype}) is just the fundamental exact sequence of
Colliot-Th\'el\`ene and Sansuc (\cite{skobook}, 
Cor.~2.3.9):
$$H^1(k,S) \to H^1(X,S) \to \Hom_k(\widehat S,\pic (\ov X)) \to H^2(k,S) 
\to H^2(X,S).$$ 

\smallskip

3. The other `extreme' case is when $\pic (\ov X)=0$. 
Then  $KD'(X)$ is quasi-isomorphic to 
$(\kbar[X]^*/\kbar^*)[1]$, and the extended type is an element of 
$\Ext^1_k(\widehat S,\kbar[X]^*/\kbar^*)$. One case of interest is
when $X$ is a principal homogeneous space of a $k$-torus $T$, so that
$\pic (\ov X)=\pic(\ov T)=0$, then the extended type is an element of
$\Ext^1_k(\widehat S,\widehat T)$. Suppose that $X=T$, and let $T'\to T$
be a surjective homomorphism of $k$-tori with kernel $S$.
This is of course a $T$-torsor under $S$. We shall show in 
Remark 2 after Proposition \ref{Pr1} below that the extended
type of this torsor is given by the natural extension
$$0 \to \widehat T \to \widehat T' \to 
\widehat S \to 0.$$ This fact was implicitly used in \cite[Lemma 2.2]{dhvoloch}.

\smallskip

4. Unlike the classical type, the extended type of a torsor 
$Y \to X$ is in general not determined by the $\ov X$-torsor $\ov Y$. 
For example, if $\pic (\ov X)=0$ and $S$ is a torus, then 
$\Ext^1_{\kbar}(\widehat S,\kbar[X]^*/\kbar^*)=0$ 
because $\widehat S$ is a free abelian group.

\smallskip

\begin{prop} \label{section}
Let $X$ be a smooth and geometrically integral variety over $k$,
and let $S$ be a $k$-group of multiplicative type.
If $X(k) \neq \emptyset$, then the map 
$\chi :  H^1(X,S) \to \Hom_k(\widehat S,KD'(X))$ is onto.
In other words, if $X(k) \neq \emptyset$, then
there exist $X$-torsors of every extended type.
\end{prop}

\dem Since $X(k) \neq  \emptyset$, the map $H^2(k,S) \to H^2(X,S)$ 
has a retraction, hence is injective. Therefore the map $\partial$ is 
zero and $\chi$ is surjective.
\enddem 

\medskip

Let $\br(X)=H^2(X,\G_{m,X})$ be the cohomological Brauer--Grothendieck 
group of $X$. As usual, $\br_0(X)$ will denote the image of 
the natural map $\br(k)\to \br(X)$, and $\br_1(X)$ the kernel of the 
natural map $\br(X)\to \br(\ov X)$.

It is easy to check that $\br_1(X)$ is canonically isomorphic
to $H^1(k,KD(X))$, see
\cite{bovh}, Prop.~2.18, or \cite{dhsz2}, Lemma~2.1. Thus
the exact triangle 
(\ref{triangle}) induces an exact sequence in Galois hypercohomology
\begin{equation}
\br(k) \to \br_1(X) \stackrel{r}{\lra} H^1(k,KD'(X)) \to H^3(k,\kbar^*).
\label{defofr}
\end{equation}
The cup-product in \'etale cohomology defines the pairing
$$\cup\ :\ H^1(k,\widehat S) \times H^1(X,S)\to
H^1(X,\widehat S)\times H^1(X,S)\to \br(X),$$
whose image visibly belongs to $\br_1(X)$. The following statement
generalises \cite[Thm. 4.1.1]{skobook}.

\begin{theo} \label{cuptheo}
Let $X$ be a smooth and geometrically integral variety over $k$,
and let $f : Y \to X$ be a torsor under a $k$-group of 
multiplicative type $S$. Let $\lambda : \widehat S \to KD'(X)$ 
be the extended type of this torsor. 
Then for any $a \in H^1(k,\widehat S)$ we have 
$$r(a \cup [Y])=\lambda_*(a),$$
where $\lambda_*$ is the induced map $H^1(k,\widehat S)\to
H^1(k,KD'(X))$.
\end{theo}

\dem 
We have canonical isomorphisms
$$H^1(X,S)=\Hom_k(\widehat S, \R p_* \G_{m,X}[1])
=\Hom_k(\widehat S, KD(X)),$$
cf. the proof of Proposition~\ref{extype}.
By \cite{milne}, Prop. V.1.20, these isomorphisms fit into the
following commutative diagram of pairings
$$\begin{array}{ccccc}
H^1(k,\widehat S)& \times& H^1(X,S)&\to&\br(X)\\
||&&||&&||\\
H^1(k,\widehat S)& \times& \Hom_k(\widehat S, \R p_* \G_{m,X}[1])&\to&
H^1(k,\R p_* \G_{m,X}[1])\\
||&&||&&\uparrow\\
H^1(k,\widehat S)& \times& \Hom_k(\widehat S, KD(X))&\to&H^1(k,KD(X))
\end{array}$$
Here the vertical arrow is induced by the canonical map
$$KD(X)=(\tau_{\leq 1}\R p_* \G_{m,X})[1]\to \R p_* \G_{m,X}[1].$$
Let $u : \widehat S \to KD(X)$ be the morphism corresponding to 
the class $[Y]$. By the commutativity of the diagram we have
$$a \cup [Y]=u_*(a) \in H^1(k,KD(X))=\br_1(X).$$
By definition, $\lambda$ is the composed map
$$\widehat S \stackrel{u}{\lra} KD(X) \stackrel{v}{\lra} KD'(X).$$
By construction,
the map $r$ from the exact sequence (\ref{defofr})
is the induced map $v_* : H^1(k,KD(X)) \to H^1(k,KD'(X))$, hence
$r(a \cup [Y])=v_*(u_*(a))=\lambda_*(a)$.
\enddem

\section{Localisation of torsors}\label{two}

Let $U$ be a smooth and geometrically integral variety over $k$.
The abelian group $\bar k[U]^*/\bar k^*$ is torsion free, 
so we can define a $k$-torus $R$ as the torus whose module
of characters $\widehat R$ is the $\Ga$-module $\bar k[U]^*/\bar k^*$.
The natural exact sequence of $\Ga$-modules
$$1\to \bar k^*\to \bar k[U]^*\to\bar k[U]^*/\bar k^*\to 1 \eqno{(*_U)}$$
defines a class 
$$[*_U]\in \Ext^1_k(\bar k[U]^*/\bar k^*,\bar k^*)=H^1(k,R).$$
Let $Y$ be the $k$-torsor under $R$ whose class in 
$H^1(k,R)$ is $-[*_U]$. 

\begin{lem} \label{le1}
There exists a morphism $q_U:U\to Y$ such that 
$q_U^*$ identifies $(*_Y)$ 
with $(*_U)$. Any morphism from $U$ to a $k$-torsor under a 
torus factors through $q_U$.
\end{lem}
\dem This is Lemma 2.4.4 of \cite{skobook}. 
\enddem

To deal with the case of torsors under arbitrary
groups of multiplicative type we need to extend
this constrution to certain geometrically reducible
varieties. However, 
Lemma \ref{le1} does not readily generalise,
because its essential ingredient is
Rosenlicht's lemma which is valid only for connected
groups. If $S$ is a torus, it says that the natural map
of $\Ga$-modules $\widehat S\to \bar k[S]^*$ induces
an isomorphism $\widehat S\cong \bar k[S]^*/\bar k^*$
(every invertible regular function that takes value $1$
at the neutral element of $S$ is a character).
This is no longer true if $\widehat S$ has
non-zero torsion subgroup.
For general groups of multiplicative type we propose the following 
substitute. 

\begin{defi} {\rm
For a (not necessarily integral) $k$-variety $V$ with an action of $S$ 
we define $\bar k[V]^*_S$ as the subgroup of 
$\bar k[V]^*$ consisting of the functions $f(x)$
for which there exists a character $\chi\in\widehat S$
such that $f(sx)=\chi(s)f(x)$ for any $s\in S(\bar k)$.}
\end{defi}

It is easy to see that $\bar k[V]^*_S$ is a 
$\Ga$-submodule of $\bar k[V]^*$.

\rem
If $S$ is a torus and $V$ is geometrically connected, 
then $\bar k[V]^*_S=\bar k[V]^*$, so we are not
getting anything new. Indeed, if $x\in V(\kbar)$
and $f\in \bar k[V]^*$, then $f(sx)/f(x)$ is a regular invertible
function on $\ov S$ with value $1$ at the neutral element
$e\in S(\kbar)$. By Rosenlicht's lemma such a function 
is a character in $\widehat S$. We obtain a morphism
from a connected variety $\ov V$ to a discrete group $\widehat S$,
which must be a constant map.
Hence there exists $\chi\in\widehat S$ such that $f(sx)=\chi(s)f(x)$
for any $x\in V(\kbar)$ and any $s\in S(\kbar)$.

\begin{prop} \label{P1}
Let $S$ be a $k$-group of multiplicative type. 
Then the natural map $\widehat S\to \bar k[S]^*_S$
induces an isomorphism of $\Ga$-modules
$\widehat S\tilde\longrightarrow\bar k[S]^*_S/\bar k^*$.
\end{prop}
\dem
The image of the natural inclusion
$\widehat S\to \bar k[S]^*$ is contained in
$\bar k[S]^*_S$, so it remains to show that any function from
$f(x)\in\bar k[S]^*_S$ that takes value $1$ at the neutral element 
$e$ of $S$
is a character. Indeed, for any $s\in S(\bar k)$ we have 
$f(sx)=\chi(s)f(x)$, and taking $x=e$ we obtain
$f(s)=\chi(s)$.
\enddem

\begin{cor} \label{c1}
Let $V$ be a $k$-torsor of $S$. Then 
we have an exact sequence of $\Ga$-modules
\begin{equation}
0\to \bar k^*\to \bar k[V]^*_S \to \widehat S\to 0. \label{E1} 
\end{equation}
The class of extension $(\ref{E1})$ in 
$\Ext^1_k(\widehat S, \bar k^*)=H^1(k,S)$ is $-[V]$.
\end{cor}
\dem
The action of $S(\bar k)$ on $\widehat S=\bar k[S]^*_S/\bar k^*$ is trivial,
hence the first statement follows from Proposition 
\ref{P1} by Galois descent. In the case when $S$ is a torus
the last statement is a well known lemma of Sansuc
\cite{sansuc}, (6.7.3), (6.7.4), see also Lemma 5.4 of \cite{bovh}.
The same calculation works in the general case.
\enddem

We shall need a relative version of Corollary \ref{c1}.
Recall that $\pi_*{\bf G}_{m,Y}$ is the sheaf on $X$ such that
for an \'etale morphism $U\to X$ we have
$$\pi_*{\bf G}_{m,Y}(U)={\rm Mor}_U(Y_U,{\bf G}_{m,U})=
{\rm Mor}_k(Y_U,{\bf G}_{m,k}).$$
Define $(\pi_*{\bf G}_{m,Y})_S$ as the subsheaf of 
$\pi_*{\bf G}_{m,Y}$ such that for an \'etale morphism $U\to X$
the group of sections  
$(\pi_*{\bf G}_{m,Y})_S(U)$ consists of the functions
$f(x)\in {\rm Mor}_U(Y_U,{\bf G}_{m,U})$ for which there exists
a group scheme homomorphism $\chi:S_U\to {\bf G}_{m,U}$ such that
$f(sx)=\chi(s)f(x)$ for any $s\in S_U(\bar k)$ and any $x\in Y_U(\bar k)$.
If $m:S_U\times_U Y_U=S\times_k Y_U\to Y_U$ is the action of $S$ on $Y_U$,
then the last condition is $m^*f=\chi\cdot f$.

\begin{prop} \label{Pr1}
Let $p:X\to{\rm Spec}(k)$ be a smooth and geometrically 
integral variety, and let $\pi:Y\to X$ be a torsor under $S$. 

{\rm (i)} 
We have an exact sequence of \'etale sheaves on $X$:
\begin{equation}
0\to {\bf G}_{m, X}  \to (\pi_*{\bf G}_{m,Y})_S \to p^*\widehat S\to 0. \label{E2} 
\end{equation}
Applying $p_*$ to $(\ref{E2})$ we obtain an exact sequence 
of $\Ga$-modules
\begin{equation}
0\to \bar k[X]^*\to \bar k[Y]^*_S\to \widehat S\to \Pic\ov X.\label{key}
\end{equation}

{\rm (ii)} The class of extension $(\ref{E2})$ in 
$\Ext^1_X(p^*\widehat S,{\bf G}_{m, X})=H^1(X,S)$ is the class 
$[Y/X]$ of the $X$-torsor $Y$ (up to sign).

{\rm (iii)} 
The last arrow in $(\ref{key})$ is the type of 
the torsor $\pi:Y\to X$. 

{\rm (iv}) When the type of $\pi:Y\to X$ is zero, the extension
given by the first three non-zero terms of $(\ref{key})$ 
maps to the class of $(\ref{E2})$ by the canonical injective map
$$0\to\Ext_k^1(\widehat S,\bar k[X]^*)\to \Ext^1_X(p^*\widehat S,{\bf G}_{m,X})=H^1(X,S).$$
\end{prop}
\dem
(i) The maps in this sequence are obvious maps.
The exactness can be checked locally, 
so we can assume that $Y=X\times_k S$, but in this case
the exactness is clear.
The exact sequence (\ref{key}) follows from (\ref{E2})
once we note that the canonical morphism $\widehat S\to p_*p^*\widehat S$
is an isomorphism since $\ov X$ is connected.

(ii) The proof of \cite[Prop. 1.4.3]{desc} applies as is.

(iii-iv) More generally, let 
$A$ be a $\Ga$-module, and ${\cal F}$ be a sheaf on $X$.
Recall that we have the spectral sequence
of the composition of functors 
$\R p_*$ and $\R \Hom_k(A,\cdot)$:
$$\Ext^m_k(A,H^n(\ov X,{\cal F}))
\Rightarrow \Ext^{m+n}_X(p^*A, {\cal F} ).$$
It gives rise to the exact sequence 
\begin{equation}
0\to \Ext_k^1(A,p_*{\cal F})\to \Ext^1_X(p^*A, {\cal F} )\to 
\Hom_k(A,R^1 p_*{\cal F}) .\label{S1}
\end{equation}
The arrows in (\ref{S1}) have explicit description.
The canonical map ${\rm E}^1\to{\rm E}^{0,1}$ sends the class of 
the extension of sheaves on $X$
$$0\to {\cal F} \to {\cal E} \to p^*A\to 0$$
to the last arrow in
$$0\to p_*{\cal F} \to p_* {\cal E} \to p_* p^*A\to R^1 p_* {\cal F},$$
composed with the canonical map $A\to p_* p^*A$. If the class of
the extension ${\cal E}$ goes to $0\in {\rm E}^{0,1}$,
then this class comes from the extension of $\Ga$-modules
$$0\to p_* {\cal F}\to p_* {\cal E} \to 
{\rm Ker}[p_*p^*A\to R^1 p_* {\cal F}]\to 0$$
pulled back by the same canonical map. See Appendix B to this paper
for a proof of these facts. In our case take $A=\widehat S$
and ${\cal F}=\G_{m,X}$. \enddem

\rems 1. The type of the torsor $\pi:Y\to X$, at least up to sign,
can also be described explicitly as follows.
Let $K=\bar k(X)$. The fibre of $p:Y\to X$ over ${\rm Spec}(K)$
is a $K$-torsor $Y_K$ under $S$. By Corollary \ref{c1}
we can lift any character $\chi\in\widehat S$ to a rational
function $f\in K[Y_K]^*_S\subset\bar k(Y)^*$. By construction $f$ is an
invertible regular function on $Y_K$, hence ${\rm div}_{\ov Y}(f)=\pi^*(D)$
where $D$ is a divisor on $\ov X$. Note that 
$D$ is uniquely determined by $\chi$ 
up to a principal divisor on $\ov X$. It is not hard to check
that the class of this divisor in $\Pic\ov X$ is the image of 
$\chi$ (up to sign). Indeed, by \cite{skobook}, Lemma 2.3.1 (ii),
the type associates to $\chi$ the subsheaf ${\cal O}_\chi$
of $\chi$-semiinvariants of $p_*({\cal O}_Y)$. The function $f$
is a rational section of ${\cal O}_\chi$, hence the class $[D]$
represents ${\cal O}_\chi\in\Pic\ov X$.
If this description was used as a definition of type, then
the exactness of (\ref{key}) is easily checked directly.

2. {From} (\ref{key}) we obtain the following exact sequence:
\begin{equation}
0\to \bar k[X]^*/\bar k^*\to \bar k[Y]^*_S/\bar k^*\to \widehat S \to \Pic\ov X.\label{key1}
\end{equation}
In the same way as in Proposition \ref{Pr1} (iii) one shows that
when the type of the torsor $\pi:Y\to X$ is zero, the extension
given by the first three non-zero terms of $(\ref{key1})$ 
maps to the extended type of $\pi:Y\to X$ by the canonical injective map
$$0\to\Ext_k^1(\widehat S,\bar k[X]^*/\bar k^*)\to \Hom_k(\widehat S,KD'(X)).$$
In particular, a surjective homomorphism of $k$-tori $T_1\to T_2$
with kernel $S$ is a $T_2$-torsor under $S$. The extended type of this torsor
comes from the extension
$$0\to\bar k[T_2]^*/\bar k^*\to\bar k[T_1]^*_S/\bar k^*\to\widehat S\to 0,$$
which is precisely the dual exact sequence
$$0\to\widehat T_2\to\widehat T_1\to \widehat S\to 0.$$

\begin{theo} \label{th1}
Let $U$ be a smooth and geometrically integral variety over $k$, 
let $S$ be a $k$-group of multiplicative type, and
let $\pi:\T\to U$ be a torsor under $S$ of type zero.
Then we have the following statements.

{\rm (i)} There is a natural exact sequence of $\Ga$-modules
$$0\to \bar k^*\to \bar k[\T]^*_S \to \widehat M\to 0,$$
which is the definition of the $k$-group of multiplicative type $M$.

{\rm (ii)} There is a natural exact sequence of $\Ga$-modules
\begin{equation}
0\to \bar k[U]^*\to \bar k[\T]^*_S\to \widehat S\to 0.\label{keyU}
\end{equation}
Let
$$1\to S\to M\to R\to 1$$
be the dual exact sequence of $k$-groups of multiplicative type. 

{\rm (iii)} We have $\T=Z\times_Y U$, where
$Z$ is a $k$-torsor under $M$ that represents the negative of 
the class of the extension {\rm (i)}, $Z\to Y=Z/S$
is the natural quotient, and
$U\to Y$ is the morphism $q_U$ from Lemma $\ref{le1}$.
\end{theo}
\dem
We note that the abelian group $\bar k[\T]^*_S/\bar k^*$
is finitely generated since the same is true for 
$\bar k[U]^*/\bar k^*$ and $\widehat S $. Thus we can define $M$
as in (i). The extension (\ref{keyU}) gives rise to
the following commutative diagram:
\begin{equation} \label{D1}
\begin{array}{ccccccccc}
&&0&&0&&&&\\
&&\downarrow&&\downarrow&&&&\\
&&\bar k^*&=&\bar k^*&&&&\\
&&\downarrow&&\downarrow&&&&\\
0&\to& \bar k[U]^*&\to& \bar k[\T]^*_S &\to &\widehat S&\to& 0\\
&&\downarrow&&\downarrow&&||&&\\
0&\to &\widehat R&\to& \widehat M&\to &\widehat S&\to &0\\
&&\downarrow&&\downarrow&&&&\\
&&0&&0&&&&
\end{array}
\end{equation}
Similarly to Lemma \ref{le1} the extension (i) defines 
a $k$-torsor $Z$ under $M$ and a morphism
$q:\T\to Z$ which identifies (i) with the extension
$$0\to \bar k^*\to \bar k[Z]^*_S\to \widehat M\to 0.$$
The functoriality of this construction and
the commutativity of (\ref{D1}) imply that 
there is an isomorphism $Z/S\cong Y$
of torsors under $R$ which makes the diagram commute
$$\begin{array}{ccc}\T&\to&U\\
\downarrow&&\downarrow\\
Z&\to&Y\end{array}$$
This gives a morphism $\T\to Z\times_Y U$ of 
$U$-torsors under $S$, which,
as any such morphism, is an isomorphism.
\enddem

\begin{cor} \label{co1}
Let $X$ be a smooth geometrically integral variety over $k$, 
let $S$ be a $k$-group of multiplicative type, 
and let $\lambda\in \Hom_k(\widehat S,KD'(X))$.
Let $U$ be a dense open set of $X$ such that 
the induced element $\lambda_U\in \Hom_k(\widehat S,KD'(U))$
has trivial image in $\Hom_k(\widehat S,\Pic \ov U)$,
so that 
$$\lambda_U\in \Ext^1_k(\widehat S,\bar k[U]^*/\bar k^*)=
\Ext^1_k(\widehat S,\widehat R)=\Ext^1_{k-{\rm groups}}(R,S),$$
where $\widehat R=\bar k[U]^*/\bar k^*$. Let 
$$1\to S\to M\to R\to 1$$
be an extension representing this class.
Then we have the following statements.

{\rm (i)} The restriction of an $X$-torsor of
extended type $\lambda$ to $U$ is isomorphic
to $Z\times_Y U$, where $Y$ is a $k$-torsor under
$R$, $U\to Y$ is the morphism $q_U$ defined in Lemma $\ref{le1}$, 
and $Z$ is a $k$-torsor under $M$ such that $Y=Z/S$.

{\rm (ii)} Conversely, any $U$-torsor $Z\times_Y U\to U$
extends to an $X$-torsor under $S$ of extended type $\lambda$.
\end{cor}
\dem
By Remark 2 after Proposition \ref{Pr1} we know that
the extension (\ref{keyU}) represents the class $\lambda_U$.
Now part (i) follows from Theorem \ref{th1}.

Recall that the embedding $j:U\to X$ gives 
a natural injective map ${\bf G}_{m,X}\to  j_*{\bf G}_{m,U}$
of \'etale sheaves on $X$. On applying 
$\R p_*$ and the truncation $\tau_{\leq 1}$ we obtain
a natural morphism 
$\tau_{\leq 1}\R p_*{\bf G}_{m,X}\to\tau_{\leq 1}\R(pj)_*{\bf G}_{m,U}$ 
in $\D(k)$. It is clear
that we have a commutative diagram of exact triangles in $\D(k)$
$$
\begin{array}{ccccc}
\bar k^* & \to & \tau_{\leq 1} \R p_*{\bf G}_{m,X} & \to &
KD'(X) [-1]\\
|| & & \downarrow & &\downarrow  \\
\bar k^* & \to & \tau_{\leq 1} \R (pj)_*{\bf G}_{m,U} & \to &
KD'(U) [-1] \end{array}
$$
It gives rise to the
following commutative diagrams of abelian groups:
\begin{equation}\begin{array}{ccccccc}
H^1(k, S)&\to & H^1(X,S) &\to&\Hom_k(\widehat S,KD'(X))&
\to & H^2(k, S) \\
||&&\downarrow&&\downarrow&&||\\
H^1(k, S)&\to&H^1(U,S) &\to&\Hom_k(\widehat S,KD'(U))&
\to &H^2(k, S)\end{array}
\label{e3} \end{equation}
Now it is easy to complete the proof of the corollary. From
Corollary \ref{c1} and Remark 2 after Proposition \ref{Pr1} 
we see that the extended type of $Z\to Y$ is $\lambda_U$.
This implies that the extended type of $Z\times_Y U\to U$ is 
$\lambda_U$. Now (ii) is an immediate consequence of (\ref{e3}).
\enddem

\section{Descent theory}
\label{three}

In this and the next chapters 
$k$ is a number field with the ring of integers $\O_k$. 
Let $\Omega_k$ be the set of places of $k$,
and let $\Omega_{\infty}$ (resp. $\Omega_f$) be the set of 
archimedean (resp. finite) places of $k$. 
For $v\in\Omega_k$ we write $k_v$ for the completion of $k$ at $v$.

For a variety $X$ over $k$ we
denote by $X(\A_k)$ the topological space of adelic points of 
$X$; it coincides with $\prod_{v \in \Omega_k} X(k_v)$ when $X$
is proper. Recall (cf. \cite{skobook}, Ch. 5) that
the Brauer--Manin pairing
$$X(\A_k) \times \br(X) \to \Q/\Z$$
is defined by the formula $$((P_v), \alpha) \mapsto 
\sum_{v \in \Omega_k} j_v(\alpha(P_v)),$$
where $j_v : \br(k_v) \to \Q/\Z$ is the local invariant in class field 
theory. By global class field theory we have 
$((P_v), \alpha)=0$ for every $\alpha \in \br_0(X)$.
For a subgroup $B\subset\br(X)$ (or $B\subset\br(X)/\br_0(X)$) 
we denote by $X(\A_k)^B$ 
the set of those adelic points that are orthogonal to $B$, 
and we write $X(\A_k)^{\br}$ for $X(\A_k)^{\br(X)}$. By the reciprocity 
law in global class field theory, we have $X(k) \subset X(\A_k)^{\br}$.

\smallskip 

If $f : Y \to X$ is a torsor under a $k$-group of multiplicative type
$S$, the {\it descent set} $X(\A_k)^f$ is defined as the 
set of adelic points $(P_v) \in X(\A_k)$ such that the family
$([Y](P_v))$ is in the image of the diagonal map
$H^1(k,G) \to \prod_{v \in \Omega_k} H^1(k_v,G)$, 
see \cite{skobook}, Section 5.3.

\begin{prop} \label{cupcor}
Let $X$ be a smooth and geometrically integral variety over 
a number field $k$,
and let $S$ be a $k$-group of multiplicative type.
An adelic point $(P_v) \in X(\A_k)$ belongs to the descent set 
$X(\A_k)^f$ associated to the torsor $f : Y \to X$ under $S$
if and only if $(P_v)$ is orthogonal to the subgroup 
$$\br_{\lambda}(X):=r^{-1}( \lambda_*(H^1(k,\widehat S)))
\subset \br_1(X)$$
with respect to the Brauer--Manin pairing.
\end{prop}

\dem The property $(P_v) \in X(\A_k)^f$ means that the family 
$([Y](P_v))$ is in the image of the diagonal map 
$H^1(k,S) \to \P^1(S)$, where $\P^1(S)$ is the restricted product 
of the groups $H^1(k_v,S)$. By the Poitou--Tate exact sequence 
(see, for example, \cite{dempoitou}, Thm. 6.3)
this is equivalent to the condition
$$\sum_{v \in \Omega_k} j_v((a \cup [Y])(P_v))=0$$
for every $a \in H^1(k,\widehat S)$. On the other hand, 
by Theorem~\ref{cuptheo} and exact sequence (\ref{defofr}) 
every element of $\br_{\lambda}(X)$ can be written as
$a \cup [Y]+\alpha_0$, 
where $\alpha_0\in \br_0(X)$. The proposition follows.
\enddem

What we want now is an `integral version' of Proposition \ref{cupcor}.
If $\Si$ is a finite set of places of $k$, we denote by $\O_\Sigma$
the subring of $k$ consisting of the elements integral at the
non-archimedean places outside $\Sigma$. Then
$U=\spec(\O_\Si)$ is an open subset of $\spec(\O_k)$.
Let us assume that there are
\begin{itemize}
\item a faithfully flat and separated $U$-scheme of finite type $\X$,
\item a flat commutative group $U$-scheme $\S$ of finite type, and
\item an fppf $\X$-torsor $\Y$ under $\S$,
\end{itemize}
such that $X=\X\times_U k$, $S=\S\times_U k$, 
and $Y=\Y\times_U k$. This assumption can always be satisfied
if $\Si$ is large enough.

Let $[\Y]$ be the class of $\Y$ in 
the fppf cohomology group $H^1(\O_{\Si},\S)=H^1(U,\S)$. 
For every $\O_{\Si}$-torsor 
$c$ under $\S$, one defines 
the {\it twisted torsor} $\Y^c=(\Y\times_U c)/\S$. This is a 
$U$-torsor under ${\cal S}$ such that
$[\Y^c]=[\Y]-c$ (see \cite{skobook}, Lemma~2.2.3).

\begin{cor} \label{cupcorentier}
Let $(P_v) \in \prod_{v \in \Si} X(k_v) \times
\prod_{v \not \in \Si} {\cal X}(\O_v)$. Then the
following conditions are equivalent: 

\smallskip

{\rm (a)} The adelic point $(P_v)$ is orthogonal to $\br_{\lambda}(X)$. 

\smallskip

{\rm (b)} 
There exists a class $[c] \in H^1(\O_{\Si}, {\cal S})$ 
that goes to
$([{\cal Y}](P_v))$ under the diagonal map
$$H^1(\O_{\Si}, {\cal S})\ \to \ \prod_{v \in \Si} H^1(k_v,S) \times 
\prod_{v \not \in \Si} H^1(\O_v, {\cal S}).$$

\smallskip

{\rm (c)} There exists 
an $\O_{\Si}$-torsor $c$ under ${\cal S}$ such that 
the adelic point $(P_v)$ lifts to an adelic point in
$\prod_{v \in \Si} Y^c(k_v) \times \prod_{v \not \in \Si} 
{\cal Y}^c(\O_v)$,
where $Y^c=\Y^c\times_{\O_\Si}k$ is the generic fibre of 
the twisted torsor $\Y^c$.

\end{cor}

\dem Let $c$ be an $\O_{\Si}$-torsor under ${\cal S}$ with 
cohomology class $[c] \in H^1(\O_{\Si},\S)$. Then $[c]$ goes to
$([\Y(P_v)])$ if and only if
$[\Y^c(P_v)]=0$ for every place $v$. 
But this is equivalent to the fact that $P_v$ lifts 
to a point in $Y^c(k_v)$ if $v \in \Si$, and to a point in $\Y^c(\O_v)$
if $v \not \in \Si$. 
This proves the equivalence of (b) and (c). 

\smallskip

Condition (b) implies that the adelic point $(P_v)$ 
is in the descent set $X(\A_k)^f$. Hence (a) follows from (b) by Proposition~\ref{cupcor}.

\smallskip 

Assume condition (a). By Proposition~\ref{cupcor}, the element 
$$([\Y(P_v)])\ \in \
\prod_{v \in \Si} H^1(k_v,S) \times \prod_{v \not \in \Si} 
H^1(\O_v, {\cal S})$$
is in the diagonal image of some $\sigma \in H^1(k,S)$.
Since $\sigma$ is 
unramified outside $\Si$, Harder's lemma (\cite{harder}, Lemma~4.1.3 or 
\cite{pianzgille}, Corollary A.8) 
implies that $\sigma$ is in the image of the restriction map 
$H^1(\O_{\Si},{\cal S}) \to H^1(k,S)$. Thus (a) implies (b).
\enddem

\rems 1. If we assume further that ${\cal S}$ and ${\cal X}$  
are smooth over $U$, then everywhere in the previous 
corollary we can replace fppf cohomology 
by \'etale cohomology. This can be arranged by choosing a
sufficiently large set $\Si$. 

2. We refer the reader to \cite{KT} and 
\cite{ctwitt} for examples of descent on
the torsor ${\cal Y}\to {\cal X}$ under $\mu_d$,
where ${\cal X}\subset{\bf P}_{\Z} ^2$ is the complement
to the closed subscheme given by
a homogeneous polynomial $f(x,y,z)$
of degree $d$ with integral coefficients,
and ${\cal Y}$ is given by the equation $u^d=f(x,y,z)$.

\medskip

Below is a ``truncated" variant of Proposition~\ref{cupcor} where we
consider all places of $k$ except finitely many.
Keep the notation as above and
let $\Si_0$ be a finite set of places of $k$.
Let $X(\A_k ^{\Si_0})$ be the topological space of ``truncated" adelic points,
defined as the restricted product of the spaces $X(k_v)$ 
for $v\not\in\Si_0$ with respect to the subsets $\X(\O_v)$, 
$v\notin \Si\cup\Si_0$. We define $\P^1_{\Si_0}(S)$ as
the restricted product of the
groups $H^1(k_v,S)$ for $v\not\in\Si_0$ with respect to
the subgroups $H^1(\O_v,\S)$, $v\notin \Si\cup\Si_0$. 
As in the classical case $\Si_0=\emptyset$, the sets 
$\P^1_{\Si_0}(S)$ and $X(\A_k ^{\Si_0})$ are independent of the 
choices of models $\S$ and $\X$. Let
$H^1_{\Si_0}(k,\widehat S)$ be the kernel of the restriction map 
$H^1(k,\widehat S)\to\prod_{v \in \Si_0}H^1(k_v,\widehat S)$.

\begin{prop} \label{var}
Let $(P_v)_{v \not \in \Si_0}\in X(\A_k ^{\Si_0})$.
Then 
$([Y](P_v))_{v \not \in \Si_0}$ is in the image of 
the diagonal map $H^1(k,S)\to\P^1_{\Si_0}(S)$ if and only if 
the ``truncated" adelic point
$(P_v)_{v\not\in\Si_0}$ is orthogonal to the subgroup 
$$\br_{\lambda,\Si_0}(X):=r^{-1}(\lambda_*(H^1_{\Si_0}(k,\widehat S)))
\subset \br_\lambda(X).$$
\end{prop}
\dem
Using the local Tate duality, we see that
the Poitou--Tate exact sequence for $S$ (see \cite{dempoitou}, Thm. 6.3.)
gives rise to the $\Si_0$-truncated exact sequence
$$H^1(k,S)\to\P^1_{\Si_0}(S)\to H^1_{\Si_0}(k,\widehat S)^{\rm D},$$
where the superscript D denotes the Pontryagin dual
$\Hom(\cdot,\Q/\Z)$. By this sequence,
$(s_v)_{v \not \in \Si_0}\in\P^1_{\Si_0} (S)$ 
is in the image of $H^1(k,S)$ if and only if
$$\sum_{v \not \in \Si_0} j_v(a \cup s_v)=0$$
for every $a \in H^1_{\Si_0}(k,\widehat S)$. 
The proof finishes in the same way as the 
proof of Proposition \ref{cupcor}. \enddem

\medskip

Taking $\Sigma=\Sigma_0$ we obtain a ``truncated" analogue of Corollary~\ref{cupcorentier}.

\begin{cor} \label{varcor}
Let $(P_v) \in \prod_{v \not \in \Si} {\cal X}(\O_v)$.
Then the following conditions are equivalent.

\smallskip

{\rm (a)} $(P_v)$ is orthogonal to $\br_{\lambda,\Si}(X)$. 

\smallskip

{\rm (b)} 
There exists a class $[c] \in H^1(\O_{\Si}, {\cal S})$ 
that goes to
$([{\cal Y}](P_v))$ under the diagonal map
$$H^1(\O_{\Si}, {\cal S})\ \to \ \prod_{v \not \in \Si} H^1(\O_v, {\cal S}).$$

\smallskip

{\rm (c)} There exists 
an $\O_{\Si}$-torsor $c$ under ${\cal S}$ such that 
$P_v$ lifts to a point in
${\cal Y}^c(\O_v)$ for every $v \not \in \Si$.

\end{cor}

Thm. 1 of \cite{dhvoloch} is a particular case of this result.

\bigskip

Let $X$ be a smooth and geometrically integral $k$-variety. Define
$$\Be(X):=\ker [\br_1 (X)/\br (k) \to \prod_{v \in \Omega_k}
\br_1 (X_v)/\br(k_v)],$$
where $X_v:=X \times_k k_v$. 
For $\alpha \in \Be(X)$ and $(P_v) \in X(\A_k)$
the image $\alpha_v$ of $\alpha$ in 
$\br_1(X_v)$ is constant for every place $v$, hence
$$i(\alpha)=\sum_{v \in \Omega_k} j_v(\alpha(P_v)) \ \in \ \Q/\Z$$ 
is well defined and does not depend 
on the choice of $(P_v)$. Let us assume that $X(\A_k) \neq \emptyset$.
Then we obtain a map $i : \Be(X) \to \Q/\Z$. Note also that this assumption,
by global class field theory,
implies that the natural map $\br(k)\to \br(X)$ is injective.
For a number field $k$ we have $H^3(k,\kbar^*)=0$, so we see from
(\ref{defofr}) that
the map $r : \br_1(X) \to H^1(k,KD'(X)) $ 
induces an isomorphism  $$\br_1(X) /\br(k) \tilde\lra H^1(k,KD'(X)).$$ 

If ${\cal C}$ is an object of $\D(k)$, and $i >0$ we define
$$\Sha^i({\cal C})=\ker [H^i(k,{\cal C}) \to 
\prod_{v \in \Omega_k} H^i(k_v,{\cal C})].$$
Thus we get an isomorphism $\Be(X)\tilde\lra\Sha^1(KD'(X))$, using
which we obtain a map $i:\Sha^1(KD'(X))\to \Q/\Z$.

\smallskip 

Let $S$ be a $k$-group of multiplicative type. There is 
a perfect Poitou--Tate pairing of finite groups
(cf. \cite{dempoitou}, Thm. 5.7)
$$\langle , \rangle_{PT} : \Sha^2(S) \times \Sha^1(\widehat S) 
\to \Q/\Z$$
defined as follows. Let
$a \in \Sha^1(\widehat S)$ and $b \in \Sha^2(S)$. 
By \cite{milne}, Lemma~III.1.16, 
the group $H^2(k,S)$ is the direct limit of the 
groups $H^2(U,{\cal S})$ where $U$ runs over non-empty open subsets 
of $\spec({\mathcal O}_k)$. Note that by taking a smaller
$U$ we can assume that $S$ and $\widehat S$ extend to 
smooth $U$-group schemes ${\cal S}$ and $\widehat {\cal S}$,
respectively. For $U$
sufficiently small we can lift $b$ to some $b_U \in H^2(U,\S)$,
and lift $a$ to some $\tilde a_U\in H^1(U,\widehat\S)$.
For any object ${\cal C}$ of $\D(U)$ we have the
hypercohomology groups with {\it compact support}
$H^i_c(U,{\cal C})$, see Section 3 of \cite{dhsza} for definitions. 
By \cite{adt}, Prop. II.2.3 (a) (see also \cite{dhsza}, Sect. 3),
since $a$ is locally trivial everywhere, $\tilde a_U$ comes from some 
$a_U \in H^1_c(U,\widehat {\cal S})$ under the natural map
$$H^1_c(U,\widehat \S)\to H^1(U,\widehat \S).$$
Define $\langle b , a \rangle_{PT}$ as the cup-product 
$b_U \cup a_U \in H^3_c(U,\G_{m,U}) \simeq \Q/\Z$ 
(the last isomorphism comes from the trace map, see 
\cite{adt}, Prop. II.2.6). It is not clear to us whether
this definition of the Poitou--Tate
pairing coincides with the classical definition in terms of cocycles,
but we shall only use the fact that it leads to
a perfect pairing.

Recall that the map $\partial:\Hom_k(\widehat S, KD'(X)) \to H^2(k,S)$ 
was defined in the exact sequence (\ref{newtype}). 

\begin{theo} \label{cuppt}
Let $X$ be a smooth and geometrically integral variety over a number
field $k$ such that $X(\A_k) \neq \emptyset$. Let $S$ be a $k$-group 
of multiplicative type,
$\lambda \in \Hom_k(\widehat S, KD'(X))$ and
$a \in \Sha^1(\widehat S)$. Then $\partial(\lambda) \in \Sha^2(S)$, and 
we have
$$\langle \partial(\lambda) , a \rangle_{PT}=i(\lambda_*(a)).$$
\end{theo}

\dem 
The image of $\partial(\lambda)$ in $H^2(X,S)$ is zero because
(\ref{newtype}) is a complex. The 
assumption $X(\A_k) \neq \emptyset$ implies that 
the map $H^2(k_v,S) \to H^2(X_v,S)$ is injective for every place $v$
(cf. also Proposition~\ref{section}).
Therefore, we have $\partial(\lambda) \in \Sha^2(S)$.

\smallskip

Recall that $w:KD'(X)\to\G_{m,k}[2]$ is the natural map
defined in (\ref{triangle}) for $X/k$.
Since (\ref{newtype}) is obtained by applying the functor 
$\Hom_k(\widehat S,.)$ to (\ref{triangle}), under the canonical
isomorphism $\Hom_k(\widehat S,\G_{m,k}[2])=H^2(k,S)$ we have the equality
$w\circ\lambda=\partial(\lambda)$. Let us write $\alpha=\lambda_*(a)\in \Sha^1(KD'(X))$.

Let $U\subset\spec(\O_k)$ be a sufficiently small
non-empty open subset such that there exists a smooth $U$-scheme
$\X$ with geometrically integral fibres and the generic fibre 
$X=\X\times_U k$, and a smooth $U$-group of multiplicative type $\S$
with the generic fibre $S=\S\times_U k$.

Write $w_U \in \Hom_U(KD'(\X),\G_{m,U}[2])$ for the map in
the exact triangle (\ref{triangle}) for $\X/U$.
The passage to the generic point $\spec(k)$ of $U$ defines 
the restriction map 
$$\Hom_U(\widehat\S,KD'(\X))\to \Hom_k(\widehat S,KD'(X)),$$
Consider the exact sequence (\ref{newtype}) for 
$\X_V=\X \times_U V$ and $\S_V=\S \times_U V$, where  
$V \subset U$ is a non-empty open set, and also for
$X$ and $S$. We obtain a commutative diagram
$$
\begin{array}{ccccccc}
H^1(\X_V,\S_V)&\to&  \Hom_V(\widehat \S, KD'(\X_V))
&\to &H^2(V,\S_V)&\to&  H^2(\X_V,\S_V) \\
\downarrow&&\downarrow&&\downarrow&&\downarrow\\
H^1(X,S)&\to& \Hom_k(\widehat S, KD'(X)) &\to&
H^2(k,S) &\to& H^2(X,S) 
\end{array}
$$
Passing to the inductive limit over $V$ and using
\cite{milne}, Lemma~III.1.16, we deduce from this diagram
a canonical surjective homomorphism
$$\varinjlim_V \ \Hom_V(\widehat \S_V,
KD'(\X_V))\to \Hom_k(\widehat S, KD'(X)).$$
Thus, by shrinking $U$, if necessary, we can lift $\lambda$ to some
$\lambda_U\in \Hom_U(\widehat\S,KD'(\X))$.
Then
$$w_U\circ\lambda_U\in \Hom_U(\widehat\S,\G_{m,U}[2])=H^2(U,\S)$$
(see the proof of Proposition~\ref{extype} for the equality here)
goes to $\partial(\lambda)$ under the restriction map
to $H^2(k,S)$.

As was explained above, we can lift $a \in \Sha^1(\widehat S)$ to some 
$a_U \in H^1_c(U, \widehat \S)$. Write $\alpha_U=\lambda_{U*}(a_U)$.
Then $\alpha_U$ is sent to $\alpha$ by the natural map
$$H^1_c(U,KD'(\X))\to H^1(k,KD'(X)).$$
By the remark before Proposition~\ref{extype} 
we can use \cite{dhsz2}, Prop.~3.3, which gives
$$\begin{array}{rcl}
i(\lambda_*(a))=i(\alpha)&=&w_U \cup \alpha_U=w_{U*}(\alpha_U)=w_{U*}(\lambda_{U*}(a_U))\\
&=&(w_{U}\circ\lambda_U)_*(a_U)=(w_{U}\circ\lambda_U)\cup a_U.
\end{array}$$
The above definition of the Poitou--Tate pairing shows 
that this equals
$\langle \partial (\lambda) , a \rangle_{PT}$.
\enddem

\rem This proof avoids delicate computations 
with cocycles as in \cite{skobook}, the proof of Thm.~6.1.2,
which follows \cite{desc}, Prop. 3.3.2.

\begin{cor}
Let $X$ be a smooth and geometrically integral variety over a number
field $k$ such that $X(\A_k)^{\Be(X)} \neq \emptyset$. Then the map 
$$\chi : H^1(X,S) \to \Hom_k(\widehat S, KD'(X))$$ is surjective
(there exist $X$-torsors of every extended type). The converse is 
true when $\pic (\ov X)$ is a finitely generated abelian group.
\end{cor}

\dem Let $\lambda \in \Hom_k(\widehat S, KD'(X)) $.
Since $X(\A_k)^{\Be(X)} \neq \emptyset$, Theorem~\ref{cuppt} 
ensures that $\langle \partial(\lambda) , a \rangle_{PT}=0$ for 
every $a \in \Sha^1(\widehat S)$. The non-degeneracy of 
the Poitou--Tate 
pairing implies that $\partial(\lambda)=0$. By 
Proposition~\ref{extype} this is equivalent to
$\lambda\in{\rm Im}(\chi)$.

\smallskip 

To prove the converse it is enough to show that 
$i:\Be(X)\to\Q/\Z$ is the zero map.
The formation of $\Be(X)$ is functorial in $X$,
so there is a natural restriction map $\Be(X^c)\to \Be(X)$.
By \cite{sansuc}, formula (6.1.4), this is an isomorphism.
The map $i^c:\Be(X^c)\to\Q/\Z$ is the composition 
$$\Be(X^c)\tilde\lra \Be(X)\stackrel{i}\lra\Q/\Z,$$
so it is enough to show that $i^c$
is identically zero.

By functoriality of the exact 
sequence (\ref{newtype}) we have a commutative diagram
with exact rows
$$\begin{array}{ccccccc}
H^1(X^c,S)& \stackrel{\chi^c}\lra &
\Hom_k(\widehat S, \pic(\ov X^c))& \stackrel{\partial^c}\lra& H^2(k,S)& 
\to& H^2(X^c,S)\\
\downarrow&&\downarrow&&||&&\downarrow\\
H^1(X,S)& \stackrel{\chi}\lra &
\Hom_k(\widehat S, KD'(X))& \stackrel{\partial}\lra& H^2(k,S)& \to& H^2(X,S)
\end{array}$$
The commutativity of this diagram implies that if $\chi$ 
is surjective, so that $\partial$ is zero, then $\partial^c$
is also zero, hence $\chi^c$ is surjective. The assumption that
$\pic(\ov X)$ is finitely generated implies that $\pic(\ov X^c)$
is also finitely generated. Using
\cite{skobook}, Prop.~6.1.4, we see that
$X^c(\A_k)^{\Be(X^c)}$ is not empty, thus
$i^c$ is identically zero.
\enddem

See \cite{witt}, Thm.~3.3.1, for miscellaneous 
characterisations of the property 
$X(\A_k)^{\Be(X)} \neq \emptyset$ in terms of the so called 
elementary obstruction and the generic period.

\section{Application: existence of integral points, 
obstructions to strong approximation}\label{four}

Recall that for a finite set of places $\Sigma_0\subset\Omega_k$ 
we denote by
$\A_k^{\Sigma_0}$ the ring of $k$-ad\`eles without $v$-components for
$v\in \Sigma_0$. 
Let $X$ be a smooth and geometrically integral $k$-variety such that 
$X(\A_k) \neq \emptyset$. There exists a finite set of places $\Sigma$
containing $\Sigma_0\cup\Omega_\infty$, and a faithfully flat morphism 
$\X\to \spec(\O_\Sigma)$ such that $X=\X\times_{\O_{\Sigma}}k$. 
We shall say that $X$ satisfies 
{\it strong approximation\footnote{We adopt the convention that
a variety $X$ 
such that $X(\A_k)=\emptyset$ satisfies strong approximation outside 
$\Sigma_0$ for every $\Sigma_0$.}
outside $\Sigma_0$} 
if $X(k)$ is dense in the restricted product 
$X(\A_k ^{\Si_0})$ 
of the sets $X(k_v)$ for $v \not \in \Sigma_0$ 
with respect to the subsets $\X(\O_v)$ (defined for $v\notin\Sigma$). The
restricted product topology is called the {\it strong} topology.
Explicitly, the base of open subsets of this topology
consists of the sets
$$\prod_{v\in T} U_v\,\times\, \prod_{v\notin T} \X(\O_v),$$
where $T$ is a finite subset of $\Omega_k\setminus\Sigma_0$
such that $\Sigma\subset T$,
and $U_v$ is an open subset of $X(k_v)$ for $v\in T$.

The following theorem gives sufficient conditions for
``the Brauer--Manin obstruction to strong 
approximation outside $\Si_0$" to be the only obstruction on $X$. 

\begin{theo} \label{strongapptheo}
Let $X$ be a smooth and geometrically integral $k$-variety such that
$X(\A_k) \neq \emptyset$, and let $S$ be a $k$-group of multiplicative 
type. Let $\Sigma_0$ be a finite set of places of $k$. 
Assume that there exists an $X$-torsor $Y$ under $S$ with the 
following property:
{\rm For all $k$-torsors $c$ under $S$, the twisted torsor 
$Y^c$ has the strong approximation property outside $\Sigma_0$.}
If $(P_v) \in X(\A_k)$ is orthogonal to 
$\br_{\lambda}(X)$, then $(P_v)_{v \not \in \Sigma_0}$ 
belongs to the closure of $X(k)$ in 
$X(\A_k^{\Sigma_0})$ for the strong topology.
\end{theo}

\dem Choose a finite set of places $\Si$ containing 
$\Si_0 \cup \Omega_{\infty}$
such that $X$ is the generic fibre of a flat smooth 
$\O_{\Si}$-scheme of finite type ${\cal X}$.
We can also assume that the torsor $Y \to X$ extends to a torsor 
${\cal Y} \to {\cal X}$ under a smooth $\O_\Si$-group scheme
of multiplicative type ${\cal S}$ such that $S=\S\times_{\O_\Si}k$. 
Furthermore, we can assume that the adelic point $(P_v) \in X(\A_k)$
belongs to $\prod_{v \in \Si} X(k_v) \times \prod_{v \not \in \Si} 
{\cal X}(\O_v)$. We want to find a rational point on $X$ 
very close to $P_v$ for $v \in (\Si-\Si_0)$ and integral outside 
$\Si$. 

\smallskip
 
By Corollary~\ref{cupcorentier}, 
the property that $(P_v)$ is orthogonal to 
$\br_{\lambda}(X)$ implies that it can be lifted to 
an adelic point $(Q_v) \in \prod_{v \in \Si} Y^c(k_v) \times 
\prod_{v \not \in \Si} {\cal Y}^c(\O_v)$ on some twisted 
torsor $Y^c$. 
In particular, $Y^c(\A_k) \neq \emptyset$.
Since $Y^c$ satisfies strong approximation 
outside $\Si_0$, we can find a rational point $m \in Y^c(k)$ 
very close to $Q_v$ for $v \in (\Si-\Si_0)$ and integral 
outside $\Si$. Sending $m$ to $X$ produces a 
rational point $m' \in X(k)$ 
very close to $P_v$ for $v \in (\Si-\Si_0)$ and 
integral outside $\Si$. 
\enddem

The following corollary gives sufficient conditions for
``the Brauer--Manin obstruction to the
integral Hasse principle" to be the only obstruction.

\begin{cor} \label{integralhasse}
Let ${\cal X}$ be a faithfully flat and separated scheme of finite type
over $\O_k$ such that $X=\X\times_{\O_k}k$.
Assume that 
$Y^c$ has the strong approximation property outside $\Omega_{\infty}$ 
for every $k$-torsor $c$ under $S$. If there exists an adelic point 
$(P_v) \in \prod_{v \in \Omega_k} {\cal X}(\O_v)$ orthogonal to 
$\br_{\lambda}(X)$, then ${\cal X}(\O_k) \neq \emptyset$.
\end{cor}

\dem Theorem~\ref{strongapptheo} 
says that $(P_v)$ can be approximated by a rational 
point $m \in X(k)$ for the strong topology on $X(\A_k ^{\Omega_\infty})$. 
Since $P_v \in 
{\cal X}(\O_v)$ for $v \in \Omega_f$, this implies that 
$m \in {\cal X}(\O_k)$. 
\enddem

As an application of Theorem \ref{strongapptheo} 
we get a short proof of a result that already 
appeared in C. Demarche's thesis \cite[Remark 4.8.2]{demthesis} (see
also \cite[Thm.~4.5]{ctfei}, where there is
an additional assumption that the geometric stabiliser $\ov H$ is finite).

\begin{theo} \label{miau}
Let $G$ be a semi-simple, simply connected linear group over a 
number field $k$. Let $\Si_0$ be a finite set of places 
of $k$ such that for every almost $k$-simple factor $G_1$ of
$G$ there exists a place $v \in \Si_0$ such that 
$G_1(k_v)$ is not compact (for example, if $k$ is not totally real
we can take $\Si_0=\{ v_0 \}$, where $v_0$ is a complex place of $k$).
Let $X$ be a homogeneous space 
of $G$ such that the geometric stabiliser $\ov H$ is a $\kbar$-group 
of multiplicative type. Then for every adelic point $(P_v)_{v \in \Omega_k}$
of $X$ orthogonal to $\br _1(X)$, the point $(P_v)_{v \not \in \Si_0}$ 
is in the closure of $X(k)$ in $X(\A_k ^{\Si_0})$
for the strong topology.
\end{theo}

In other words: the Brauer--Manin obstruction to strong
approximation outside $\Si_0$ is the only one on $X$.

\dem Let us assume that $G$ acts on $X$ on the left.
Then $\ov X$ with the left action of $\ov G$
is isomorphic to $\ov G/\ov H$. 
Since $\pic (\ov G)=0$ and $\kbar[G]^*=\kbar^*$, the abelian group
$\pic (\ov X)$ is finitely generated, and $\kbar[X]^*=\kbar^*$. Now
the existence of a point $(P_v)_{v \in \Omega_k}$ orthogonal to 
$\br_1(X)$ implies that $X(k) \neq \emptyset$ by 
\cite{skobook}, Prop. 6.1.4, and \cite{dhsko}, Prop 3.7 (3) and 
Example 3.4. Therefore $X$ with the left action of $G$
is isomorphic to $X=G/H$, where $H$ 
is a $k$-group of multiplicative type. Taking $Y=G$, we obtain
a right torsor $Y \to X$ under $H$ such that 
for any $k$-torsor $c$ under $H$ the twist
$Y^c$ is a left $k$-torsor under $G$. 

By the Hasse principle for semi-simple simply connected groups 
(a theorem of Kneser--Harder--Chernousov), $Y^c({\bf A}_k) \neq \emptyset$
implies $Y^c(k)\not=\emptyset$, hence $Y^c\simeq G$. 
By the strong approximation 
theorem (see, for example, \cite{platonov}, Thm. 7.12), $G$
satisfies strong approximation outside $\Si_0$.
It remains to apply Theorem~\ref{strongapptheo}.
\enddem

\rems 1. It is not clear to us whether Corollary~\ref{integralhasse} 
still holds if we only assume 
that all the twists $Y^c$ satisfy the integral Hasse principle: indeed, 
we do not know in general whether the torsor $Y \to X$ can be 
extended to an fppf torsor ${\cal Y} \to {\cal X}$ over 
$\spec (\O_k)$. 

\smallskip

2. The assumptions of 
Theorem~\ref{strongapptheo} and Corollary~\ref{integralhasse} imply that
$\kbar[X]^*=\kbar^*$, that is, we are still in ``the classical case" 
of descent theory. Indeed, if $Y$ satisfies 
strong approximation outside a finite set of places, then 
$\ov Y$ is simply connected (this was first observed in
\cite{minchev}, Thm. 1, see 
also \cite{dhens}, Cor.~2.4). This implies $\kbar[Y]^*=\kbar^*$,
and hence $\kbar[X]^*=\kbar^*$.
(Otherwise pick up a function
$f\in \kbar[Y]^*$ such that the image of $f$ in the free abelian group
$\kbar[Y]^*/\kbar^*$ is not divisible by a prime $\ell$.
Then the normalisation of $\ov Y$ in $\kbar(Y)(f^{1/\ell})$ 
is a connected \'etale covering of $\ov Y$ of degree $\ell$.)

\section*{Appendix A}

Let $Z$ be an integral regular Noetherian scheme, and let
$p:X\to Z$ be a smooth faithfully flat morphism of finite type
with geometrically integral fibres.
The goal of this appendix is to show that the object
$\tau_{\leq 1} \R p_* \G_{m,X}$ of the derived category $\D(Z)$
of \'etale sheaves on $Z$
can be represented by an explicit two-term complex.
This links our $KD(X)$ and $KD'(X)$ with their analogues
introduced in \cite{dhsz2}, Remark 2.4 (2). 

\medskip

Let $j: \eta=\Spec (k(X)) \hookrightarrow X$ be the inclusion
of the generic point. Since $X$ is regular, there is no difference
between Weil and Cartier divisors, so we have the following exact sequence
of sheaves on $X$, see \cite{milne}, Examples II.3.9 and III.2.22:
$$
0\to \G_{m,X} \to j_*\G_{m,\eta} \to \Div_X\to 0,
$$
where $\Div_X$ is the sheaf of divisors on $X$, that is,
the sheaf associated to the presheaf such that the group of sections over
an \'etale $U/X$ is the group of divisors on $U$.

We call an irreducible effective divisor $D$
on $X$ {\it horizontal} if it is
the Zariski closure of a divisor on the generic fibre of $p:X\to Z$.
If $D=p^{-1}(D')$ for a divisor $D'$ on $Z$, we call
$D$ {\it vertical}. The sheaf
$\Div_X$ is the direct sum of sheaves
$$\Div_X=\Div_{X/Z}\oplus \Div^v_X,$$
where $\Div_{X/Z}$ is the subsheaf of horizontal divisors, and 
$\Div^v_X$ is the subsheaf of vertical divisors.

Define a subsheaf $K^\times_{X/Z}\subset j_*\G_{m,\eta}$ by the 
condition that the following diagram is commutative and has
exact rows and columns:
$$
\begin{array}{ccccccccc}
&&&&0&&0&&\\
&&&&\downarrow&&\downarrow&&\\
0&\to& \G_{m,X}& \to &K^\times_{X/Z}& \to &\Div_{X/Z}&\to& 0\\
&&||&&\downarrow&&\downarrow&&\\
0&\to& \G_{m,X}& \to &j_*\G_{m,\eta}& \to &\Div_X&\to& 0\\
&&&&\downarrow&&\downarrow&&\\
&&&&\Div^v_X&=&\Div^v_X&&\\
&&&&\downarrow&&\downarrow&&\\
&&&&0&&0&&
\end{array}
$$
The complex of \'etale sheaves on $Z$ 
$$p_*K^\times_{X/Z}\to p_*\Div_{X/Z},$$
after the shift by 1 to the left,
is the complex $\K\D(\X)$ defined in \cite{dhsz2}, Remark 2.4 (2),
see also the formulae on the bottom of page 538. There is a 
natural injective morphism $\G_{m,Z}\to p_*K^\times_{X/Z}$;
the complex
$$p_*K^\times_{X/Z}/\G_{m,Z}\to p_*\Div_{X/Z}$$
was introduced in \cite{dhsz2} and denoted there by $\K\D'(\X)$.

\medskip

\noindent{\bf Proposition\,} {\it
The object 
$\tau_{\leq 1} \R p_* \G_{m,X}$ of the derived category of
\'etale sheaves on $Z$ is represented by the complex
$p_*K^\times_{X/Z}\to p_*\Div_{X/Z}$.}
\dem
The proof of Lemma 2.3 of \cite{bovh} works in our situation.
To complete the proof we only need to show that
$R^1p_*K^\times_{X/Z}=0$. Note that the canonical morphism $\Div_Z\to p_*\Div^v_X$
is an isomorphism because $p$ is surjective with geometrically integral 
fibres. Now the exact sequence of sheaves on $X$
$$0\to K^\times_{X/Z}\to j_*\G_{m,\eta}\to \Div^v_X\to 0$$
gives rise to the following exact sequence of sheaves on $Z$:
$$p_*j_*\G_{m,\eta}\to \Div_Z\to R^1p_*(K^\times_{X/Z})\to R^1p_*(j_*\G_{m,\eta}).$$
Using the spectral sequence of the composition of functors 
$\R p_*$ and $\R j_*$ we see that the sheaf $R^1p_*(j_*\G_{m,\eta})$ 
has a canonical embedding into $R^1(pj)_*\G_{m,\eta}$.
The latter sheaf is zero by Grothendieck's version of Hilbert's theorem 90.

It remains to prove the surjectivity of $(pj)_*\G_{m,\eta}\to \Div_Z$,
which is enough to check at the stalk at any geometric point of $Z$.
But locally every divisor on $Z$ is the divisor of a function,
since $Z$ is regular. This completes the proof. \enddem

\rem In this appendix we worked over the small \'etale site of $Z$. 
Applying our arguments to an arbitrary smooth scheme of finite 
type $S/Z$ one shows that the same results remain true
for the smooth site {\it Sm}/$Z$ used in \cite{dhsz2}.

\section*{Appendix B}

The functor $\R \Hom_X(p^*A,\cdot):\D(X)\to\D({\rm Ab})$ 
is the composition of functors
$\R p_*:\D(X)\to\D(k)$ and $\R\Hom_k(A,\cdot):\D(k)\to \D({\rm Ab})$, 
hence we have
$$
\R\Hom_X(p^*A,{\cal F})=\R\Hom_k(A,\R p_*{\cal F}).
$$
Explicitly, this isomorphism associates to $p^*A \to {\cal F}$ the composition
$$A\to \R p_* (p^*A)\to \R p_*{\cal F},$$
where the first map is the canonical adjunction morphism. The inverse
associate to $A\to \R p_*{\cal F}$ the composition 
$$p^*A\to p^*(\R p_*{\cal F})\to {\cal F},$$ where the last map
is the second canonical adjunction morphism.

\medskip

Let us now complete the proof of Proposition \ref{Pr1} (iii).
To give an equivalence class of the extension of sheaves on $X$
\begin{equation}
0\to {\cal F}\to {\cal E}\to p^*A\to 0\label{f4}
\end{equation}
is the same as to give a morphism $p^*A \to {\cal F}[1]$ in 
the derived category $\D(X)$. By the above, to this morphism we
associate the composition
$$A\to \R p_* p^*A\to \R p_* {\cal F}[1].$$ Since $A$ is a one-term
complex concentrated in degree 0 this composition comes from a morphism
$\alpha:A\to(\tau_{\leq 1}\R p_*{\cal F})[1]$ in $\D(k)$.
By taking the $0$-th cohomology 
we obtain a homomorphism $\beta:A\to R^1p_*{\cal F}$ of discrete Galois
modules. Clearly, $\beta$ is the composition of 
the canonical map $A\to p_*p^*A$
with the differential in the long exact sequence
of cohomology attached to (\ref{f4}):
$$0\to p_*{\cal F}\to p_*{\cal E}\to p_*p^*A\to R^1p_*{\cal F}.$$
To finish the proof of (iii) we need to show that
$\beta$ can also be obtained through the spectral sequence, that is,
as the image of the class of (\ref{f4}) under the right arrow in (\ref{S1}).
But (\ref{S1}) is obtained by applying $\R\Hom_k(A,\cdot)$ to
the exact triangle
\begin{equation}
(p_*{\cal F})[1]\to(\tau_{\leq 1}\R p_*{\cal F})[1]\to R^1 p_*{\cal F}.\label{f2}
\end{equation}
By definition, $\beta$ is the composition of $\alpha$ 
with the right map in (\ref{f2}), so the proof of (iii) is now complete.

\medskip

Let us complete the proof of Proposition \ref{Pr1} (iv).
The exact triangle (\ref{f2}) gives rise to the exact 
sequence of abelian groups
$$0\to \Hom_k(A,(p_*{\cal F})[1])\to 
\Hom_k(A,(\tau_{\leq 1}\R p_*{\cal F})[1])
\to \Hom_k(A,R^1p_*{\cal F}),$$
which is the same as (\ref{S1}). 
Since the right arrow here sends $\alpha$ to $\beta$, we see that
if $\beta=0$, then $\alpha$ comes from
a morphism $A\to (p_*{\cal F})[1]$. Hence the class of 
(\ref{f4}) comes from
the class of an extension of $A$ by $p_*{\cal F}$, say
\begin{equation}
0\to p_*{\cal F}\to B\to A\to 0, \label{f5}
\end{equation}
in the sense that (\ref{f4}) is the push-out of 
$$0\to p^*p_*{\cal F}\to p^* B\to p^*A\to 0$$
by the adjunction map $p^*p_*{\cal F}\to {\cal F}$. Therefore,
by the description of the adjunction isomorphism and its inverse
given above, applying $p_*$ to (\ref{f4}), and pulling back the resulting
short exact sequence via the adjunction map $A\to p_*p^*A$ 
(this makes sense when $\beta=0$) gives back the extension (\ref{f5}). 
\enddem

\bigskip

\noindent {\bf Acknowledgement} We would like to thank the organisers of
the workshop on anabelian geometry at the Isaac Newton Institute 
for Mathematical Sciences in Cambridge (August, 2009), 
where the work on this paper began.

\bigskip

\noindent
Math\'ematiques, B\^atiment 425, Universit\'e Paris-Sud,
Orsay, 91405 France
\medskip

\noindent David.Harari@math.u-psud.fr

\bigskip

\bigskip

\noindent Department of Mathematics, South Kensington Campus,
Imperial College London,
SW7 2BZ England, U.K.

\smallskip

\noindent Institute for the Information Transmission Problems,
Russian Academy of Sciences, 19 Bolshoi Karetnyi,
Moscow, 127994 Russia
\medskip

\noindent a.skorobogatov@imperial.ac.uk

\begin{thebibliography}{99}

\bibitem{boduke} M. Borovoi. Abelianization of the second 
nonabelian Galois cohomology. {\it Duke Math. J.} {\bf 72} (1993) 217--239.

\bibitem{bovh} M. Borovoi and J. van Hamel. 
Extended Picard complexes 
and linear algebraic groups. {\it J. reine angew. Math.}  
{\bf 627} (2009) 53--82.

\bibitem{blr} S. Bosch, W. L\"utkebohmert and M. Raynaud. 
{\it N\'eron models}. 
Ergebnisse der Mathematik und ihrer Grenzgebiete,
Springer-Verlag, 1990.

\bibitem{breen2} L. Breen. 
On a nontrivial higher extension of representable abelian sheaves. 
{\it Bull. Amer. Math. Soc.} {\bf 75} (1969) 1249--1253.

\bibitem{breen} L. Breen. 
Un th\'eor\`eme d'annulation pour certains $Ext^{i}$ 
de faisceaux ab\'eliens.  {\it Ann. Sci. \'Ecole Norm. Sup.} 
{\bf 8} (1975) 339--352.

\bibitem{ctfei} J-L. Colliot-Th\'el\`ene and Xu Fei.
Brauer--Manin obstruction for integral points of 
homogeneous spaces and representation of integral quadratic forms.
{\it Comp. Math.} {\bf 145} (2009) 309--363.

\bibitem{desc} J-L. Colliot-Th\'el\`ene et J-J. Sansuc. 
La descente 
sur les vari\'et\'es rationnelles, II. {\it Duke Math. J.} {\bf 54} 
(1987) 375--492.

\bibitem{ctwitt} J-L. Colliot-Th\'el\`ene et O. Wittenberg.
Groupe de Brauer et points entiers de deux familles de surfaces 
cubiques affines. {\it Amer. J. Math}, to appear. 
Available at
http://www.math.ens.fr/\verb1~1wittenberg/troiscubes.pdf

\bibitem{demthesis} C. Demarche. 
{\it M\'ethodes cohomologiques pour l'\'etude 
des points rationnels sur les 
espaces homog\`enes}. Th\`ese de l'Universit\'e Paris-Sud, 2009. 
Available at 
http://www.math.u-psud.fr/\verb1~1demarche/thesedemarche.pdf

\bibitem{dempoitou} C. Demarche. Suites de Poitou--Tate pour les
complexes de tores \`a deux termes. {\it I.M.R.N.}, to appear.
Available at
http://www.math.u-psud.fr/\verb1~1demarche/dualitetores.pdf

\bibitem{dem2} C. Demarche. Le d\'efaut d'approximation forte dans 
les groupes lin\'eaires connexes.
{\it Proc. London Math. Soc.}, to appear. arXiv:0906.3456v1

\bibitem{pianzgille} P. Gille and A. Pianzola.
Isotriviality and \'etale cohomology of Laurent polynomial rings.
{\it J. Pure Appl. Algebra}  {\bf 212} (2008) 780--800.

\bibitem{sga3} A. Grothendieck et al. {\it Sch\'emas en Groupes}
(SGA 3), I, Lecture Notes in Math. {\bf 151}, Springer-Verlag, 1970.

\bibitem{dhens} D. Harari. Weak approximation and non-abelian 
fundamental groups. {\it Ann. Sci. E.N.S.} {\bf 33} (2000) 467--484.

\bibitem{dhsko} D. Harari and A. N. Skorobogatov. 
Non-abelian cohomology and rational points.
{\it Comp. Math.} {\bf 130} (2002) 241--273.

\bibitem{dhsza} D. Harari and T. Szamuely. Arithmetic duality 
theorems for 1-motives. {\it J. reine angew. Math.}  {\bf 578} (2005) 93--128.

\bibitem{dhsz2} D. Harari and T. Szamuely. Local-global principles 
for 1-motives. {\it Duke Math. J.}  {\bf 143} (2008) 531--557.

\bibitem{dhvoloch} D. Harari and J.F. Voloch. 
The Brauer--Manin obstruction for integral points on curves. 
{\it Math. Proc. Cambridge Philos. Soc.}  {\bf 149} (2010) 413--421.

\bibitem{harder} G. Harder. Halbeinfache Gruppenschemata 
\"uber Dedekindringen. {\it Inv. Math.}  {\bf 4} (1967) 165--191.

\bibitem{KT} A. Kresch and Yu. Tschinkel. Two examples of Brauer--Manin
obstruction to integral points. {\it Bull. London Math. Soc.} {\bf 40}
(2008) 995--1001.

\bibitem{adt} J.S. Milne. {\it Arithmetic duality theorems}, 
Academic Press, 1986.

\bibitem{milne} J.S. Milne. {\it \'Etale cohomology}, Princeton University
Press, 1980.

\bibitem{minchev} Kh. P. Minchev. Strong approximation for
varieties over algebraic number fields. 
{\it Dokl. Akad. Nauk BSSR} {\bf 33} (1989) 5--8. (Russian)

\bibitem{mum} D. Mumford. {\it Abelian varieties}, 
Oxford University Press, 1970.

\bibitem{platonov} V. Platonov and A. Rapinchuk. {\it Algebraic 
groups and number theory}, Academic Press, 1994.

\bibitem{sansuc} J-J. Sansuc. Groupe de Brauer et arithm\'etique 
des groupes alg\'ebriques lin\'eaires sur un corps de nombres.
{\it J. reine angew. Math.} {\bf 327} (1981) 12--80.

\bibitem{S99} A.N. Skorobogatov. Beyond the Manin obstruction.
{\it Inv. Math.} {\bf 135} (1999) 399--424.

\bibitem{skobook} A.N. Skorobogatov. {\it Torsors and rational points}, 
Cambridge Tracts in Mathematics {\bf 144}, 
Cambridge University Press, 2001.

\bibitem{witt} O. Wittenberg. On Albanese torsors and the elementary 
obstruction. {\it Math. Ann.} {\bf 340} (2008) 805--838.

\end{thebibliography}
\end{document}